\documentclass[reqno]{amsart}
\usepackage{amsmath}
\usepackage{amssymb,amsmath,amsthm,amsfonts,color,breqn,pgf}  
 \usepackage[abs]{overpic} 
\usepackage{xcolor}
\usepackage[colorlinks=true, pdfstartview=FitV, linkcolor=blue, 
citecolor=blue, urlcolor=blue]{hyperref}
\usepackage{autonum}

\makeatletter
\@namedef{subjclassname@2020}{\textup{2020} Mathematics Subject Classification}
\makeatother

\theoremstyle{plain}
\newtheorem{theorem}{Theorem}[section]
\newtheorem{proposition}[theorem]{Proposition}
\newtheorem{lemma}[theorem]{Lemma}

\newcommand\epsi{\varepsilon}
\newcommand\ove{\overline}
\newcommand\haz{\widehat}
\newcommand{\Rz}{{\mathbb R}}

\renewcommand{\d}{{\rm d}}

\newcommand\UUU{\color{black}}
\newcommand\TTT{\color{black}}

\newcommand\EEE{\color{black}}
\newcommand\RRR{\color{black}}

\title[John regularity]{Level sets of  solutions to the stationary
  \UUU Hamilton--Jacobi \EEE equation  are John regular}
\author[E. Davoli] {Elisa Davoli} 
\address[Elisa Davoli]{Institute of Analysis and Scientific Computing, TU Wien, 
Wiedner Hauptstrasse 8-10, A-1040 Vienna, Austria}
\email{elisa.davoli@tuwien.ac.at}
\author[U. Stefanelli]{Ulisse Stefanelli} 
	\address[Ulisse Stefanelli]{University of
		Vienna, Faculty of Mathematics,
                Oskar-Morgenstern-Platz 1, A-1090 Vienna, Austria, 
		University of Vienna \UUU and \EEE Vienna Research Platform on Accelerating
		Photoreaction Discovery, W\"ahringerstra\ss e 17, 1090
                Wien, Austria}
	\email{ulisse.stefanelli@univie.ac.at}
	\urladdr{http://www.mat.univie.ac.at/$\sim$stefanelli}

\subjclass[2020]{35F21, 
35B65
}
\keywords{\UUU Hamilton--Jacobi \EEE equation, viscosity solution, regularity, John domain.}

\begin{document}

\begin{abstract} Let $u$ be the unique nonnegative viscosity solution of
   the Hamilton--Jacobi equation $H(x,\UUU \nabla u \EEE )=0$ in the external domain $\Rz^{ n} \setminus K$
  with $u=0$ on $K$. Under general conditions on $H$, we prove
  that all sublevels of $u$ are John domains.
  Moreover, if $K$
  itself is a John domain, we  provide  a uniform lower bound on
  the John constant of all sublevels. We  exhibit  counterexamples showing that
  John regularity is sharp in this setting.
\end{abstract}
 
\maketitle

\section{Introduction}
This note is concerned with the regularity of solutions to
the external
problem for the stationary Hamilton--Jacobi equation, 
namely, 
  \begin{equation}
     H(x,\UUU\nabla u\EEE)= 0\quad \text{in} \ \Rz^{ n} \setminus K,
                           \quad 
    u=0\quad \text{on} \ K.\label{eq:0}
  \end{equation}
This problem \UUU arises in different contexts, \EEE from 
front propagation \cite{Sethian,Tran}, to geometric optics
\cite{Kline}, optimal control \cite{Bardi,Cannarsa2}, differential
games \cite{Yong}, image processing \cite{Aubert}, information flow
\cite{Dunbar}, and growth \cite{Caflisch}. 

We consider the general case of a 
 discontinuous Hamiltonian $H:\Rz^n \times \Rz^n \to \Rz$, which is
convex \UUU and \EEE  nondegenerate \UUU in the
gradient \EEE variable. Our setting, made precise by \eqref{eq:H} later
on, covers the reference case of the {\it generalized eikonal
  equation}
$$H(x,p)=\alpha(x)|p|-1 \quad
\text{with}\quad 0<\alpha_*\leq \alpha(\cdot)\leq
\alpha^*,$$
\UUU for $\alpha$ merely continuous, \EEE and guarantees that, for all given nonempty compact $K \subset \Rz^{ n}$, problem \eqref{eq:0}
admits a unique nonnegative viscosity solution  $u$, cf.~
\cite[Thm.~3.15]{Mennucci} and Section \ref{sec:main}. 


We \UUU investigate \EEE the regularity of the open
sublevels  
$$U_t:=\{x \in \Rz^{ n}\: :\: u(x) < t \}\quad\text{\UUU for} \  t>0.$$
Our interest is
motivated by the relevance of the sets $U_t$ in  many  different
applications \cite{Falcone}. In front propagation, $U_t$
corresponds to the points reached  at  time $t$ by a front originated
at $K$ \cite{Sethian}. In optimal control, $U_t$ is interpreted as  {\it reachable set}
for a controlled ODE system \cite{Bardi}, see also Section \ref{sec:counter}. In modeling growth \cite{Caflisch,morpho3},
$U_t$ corresponds to the region occupied at time $t$ by \UUU body
growing by accretion. \EEE 

In order to introduce our regularity result, let us first
consider the
case of the classical {\it eikonal
  equation} $|\UUU \nabla u\EEE|=1$. In this case, $u$ is simply the 
distance to $K$, namely, $u(x)=d(x,K):=\inf_{y\in
  K}|x-y|$. In particular, for all $t>0$ one has that $U_t = K +
B_{t}(0)$.  Thus, the sets $U_t$  are {\it interior-ball regular},
\UUU where we recall that a set $U\subset \Rz^n$ is interior-ball
regular if \EEE there exists a radius $r>0$ such that for all $x\in
\partial U$ one finds $y \in U$ with $|x-y|=r$  so that 
$B_r(y):=\{w\in \Rz^n\::\: |w-y|<r\}\subset U$.  
In fact,  the  
interior-ball-regularity radius of $U_t$ is exactly $t$. Note $K$ is not
required to be regular. In particular, the eikonal equation
smoothens $U_t$ instantaneously.

In the case of the generalized eikonal equation $\alpha(x)|\UUU\nabla
u\EEE|=1$ with nonconstant $\alpha$, the interior-ball
regularity of superlevels still holds, as long as one
assumes that $\alpha \in
C^{1,1}(\Rz^{ n})$. This
fact has been proved by {\sc Lorenz} \cite{Lorenz}. 
The reader is
also referred to  \cite{Alvarez,Cannarsa,Caroff} \UUU and to
\cite{Barles09,Chen} \EEE for
closely related interior-ball-regularity results  for smooth $\alpha$,
\UUU also in the parabolic case. \EEE Note
nonetheless that the  interior-ball-regularity
radius may degenerate
as $t \to \infty$, \UUU see Section \ref{sec:counter1} below. \EEE

The aim of this note is to \UUU prove sharp regularity results for
the  \EEE
sublevels $U_t$ in full generality, 
covering the case of discontinuous Hamiltonians, \UUU as well. \EEE In this
situation, the interior ball-regularity of $U_t$ cannot be expected,
as we show in
Section \UUU \ref{sec:counter3}. \EEE
Our main result, Theorem \ref{thm:main},
states that the sublevels $U_t$ are nonetheless John domains for all
$t >0 $. 

We recall that a  nonempty  bounded domain $U \subset
\Rz^{ n}$ is said to be  a  
{\it John domain} (equivalently,   {\it John regular}) with respect to a  fixed  point $x_0 \in U$ and a given
{\it John constant} $\kappa\UUU \in (0,1]$ \EEE if it satisfies an internal {\it twisted cone condition}: for all points $x\in U$ one can find an
arc-length parametrized
curve $\rho: [0,L_\rho] \to U$ such that $\rho(0)=x$,
  $\rho(L_\rho) = x_0$, and $ d  (\rho(s), \partial U) \geq
  \kappa \, s$ for all $s\in [0,L_\rho]$.
  Note that John domains are connected. 

 \UUU More precisely, \EEE we prove that, starting from a John domain, the
  sublevels $U_t$ are John domains with respect to the some uniform John
  constant, independently of $t>0$.
Our main result reads as follows. 

\begin{theorem}[John regularity]\label{thm:main} Let $H:\Rz^n
  \times \Rz^n \to \Rz $ be such that
  \begin{align}
  & C:=\{(x,p)\in \Rz^n \times \Rz^n\::\: H(x,p)\leq 0\}\ \text{is closed}, \nonumber\\[1mm]
  &C_x:=\{p\in \Rz^n\::\:
    H(x,p)\leq 0\}\ \text{is \UUU  convex} \quad \forall
    x \in \Rz^n, \nonumber \\
  &C_x= \overline{\{ p \in \Rz^n \::\: (x,p)\in \,C^\circ\}}\quad \forall
    x \in \Rz^n, \nonumber \\[1mm]
  & \exists 0<\alpha_* \leq \alpha^*: \quad B_{1/\alpha^*}(0)\subset
    C_x \subset  B_{1/\alpha_*}(0)\quad \forall x \in\Rz^n. \label{eq:H}
  \end{align}
  Let
  $K\subset \Rz^n$ be compact and connected and $u$ be the unique \UUU
  nonnegative viscosity solution to \EEE \eqref{eq:0}. Then, the sublevels $U_t$ are John domains for all
  $t>0$.

  \UUU Moroever, if \EEE  $K$ is a  John domain with respect   $x_0\in  K^\circ$ with John constant
  $\kappa_0>0$,  all
  sublevels $U_t$  are John
  domains with respect to $x_0$ with John constant 
\begin{equation}
  \kappa :=\frac{\alpha_*}{2\alpha^*+\alpha_*}\UUU \kappa_0. \EEE
\label{eq:CJohn}
\end{equation}
\UUU For $t> 2 \sup\{r>0\,:\, B_r(x_0)\in K\}/(\alpha_*\kappa_0)$, one
can choose $\kappa=\alpha_*/(2\alpha^*+\alpha_*)$. \EEE
\end{theorem}

\UUU In the statement above, we have used the symbols $\ove E$ and
$E^\circ$ to indicate the closure and the interior of a set $ E$,
respectively. In the special case of the generalized eikonal equation
$H(x,p)=\alpha(x)|p|-1$, assumption \eqref{eq:H} implies that $\alpha$
is continuous. Hence,
$H$ is continuous, as well. Note however that assumption  \eqref{eq:H}
does not imply continuity of $H$ in general: take $H(x,p) =
h(x)(|p|-1)$ for $h:\Rz^n \to (0,\infty)$ arbitrary. In this case, $C=
\{(x,p) \in \Rz^n \times \Rz^n\,:\, |p|\leq 1\}$, $C^\circ=\{(x,p) \in
\Rz^n \times \Rz^n\,:\, |p|< 1\}$, $\{p\in \Rz^n
 \,:\,  (x,p)\in C^\circ\}= B_1(0)$, and 
$C_x=\ove{B_1(0)}$. Hence, \eqref{eq:H} follows, although $H$ is not
necessarily continuous. \EEE

The proof of Theorem \ref{thm:main} is given in Section \ref{sec:main}
below. \UUU Our argument is based on a geometric construction using
the fact that the Hamilton--Jacobi equation induces a Finsler metric on $\Rz^n$ \cite{Siconolfi}
which is equivalent to the Eulerian one.\EEE

 In Section \ref{sec:counter}, we
provide some examples illustrating the sharpness  of
Theorem~\ref{thm:main}.   In particular, in Section
\ref{sec:sharpK} we show  that no uniform lower
bound on the John constant can be expected if $K$ is not John
regular.

In the very general setting of Theorem \ref{thm:main}, \UUU no \EEE
interior-ball (nor  interior-cone)
regularity of the sublevels $U_t$ can be expected in
general.  We provide explicit examples of
this fact in
Sections \ref{sec:counter3}--\ref{sec:counter4} in the special
case of the generalized eikonal equation $\alpha(x)|\nabla \UUU u \EEE
|=1$
 for a smooth $K$ and for $\alpha$ Lipschitz continuous but not 
$C^{1,1}$. As a by product, \UUU this proves \EEE that the $C^{1,1}$- 
regularity requirement for $\alpha$ in \cite{Alvarez,Cannarsa,Caroff,Lorenz} is
necessary for
interior-ball regularity.

In contrast to to the interior-ball-regularity case, \UUU note that
the \EEE  John regularity \UUU of $U_t$ \EEE does not degenerate \UUU
as $t \to \infty$. \EEE 
Our lower bound on the John constant from \eqref{eq:CJohn}, albeit expectedly
not optimal, holds for all $t>0$. \UUU Moreover, starting from some
given time, also depending on the initial set $K$, the  
John constant $\kappa$ of $U_t$ can be even chosen independently from
the John constant $\kappa_0$ of $K$. \EEE 

 Before closing this introduction, let us mention some
 consequences of our result. 
In the context of optimal control, we mention the papers
\cite{Boarotto,Duprez} focusing on the case of  a  generalized eikonal
equation and proving that the boundary $\partial U_t$ of the
reachable set $U_t$ is 
negligible with respect to the $n$-dimensional Lebesgue measure, under different regularity
requirements for $K$.
 Albeit in a slightly different setting,  our result provides an alternative, more general, and sharper take to this fact: By checking that
$U_t$ is John regular (in fact, even in the more general case of
discontinuous Hamiltonians) and by applying \cite[Cor. 2.3]{Koskela} one
directly obtains that
\begin{equation}
  {\rm dim}_{M}(\partial U_t) \leq n - \UUU \mu\EEE<n,\label{eq:koskela}
\end{equation}
where ${\rm dim}_{M}$ denotes the classical Minkowski or
box-counting dimension, and $\UUU \mu \EEE >0$. This in particular entails
that $\mathcal H^s (\partial U_t)=0$ for all $s>n-\UUU \mu \EEE$, where $\mathcal H^s$ is the
$s$-dimensional Hausdorff measure. By taking $s=n$ we recover \UUU a
stronger result with respect to \EEE \cite{Boarotto,Duprez}.

In Section \ref{sec:regu}, we revisit \eqref{eq:koskela} in order to
\UUU address \EEE  the regularity of the set-valued map $t\mapsto U_t$. In
particular, we have the following.

\begin{proposition}[Regularity of the map $t \mapsto U_t$]\label{prop:regu}
Under assumption \eqref{eq:H}, let
  $K\subset \Rz^n$ be compact and connected, and let $u$ be the unique
  \UUU nonnegative viscosity solution to \EEE ~\eqref{eq:0}. Then,
  \begin{align}
    &\alpha_*(t-r)\leq \inf_{x\in \partial U_t}d(x,\partial
    U_r) \quad \forall 0<r\leq t. \label{eq:lipschitz0}\\
    &\sup_{x\in U_t}d(x,U_r)\leq \alpha^* (t-r) \quad \forall 0<r\leq t.  \label{eq:lipschitz}
  \end{align} 

  If $K$ is a  John domain with respect to  \UUU $x_0\in  K^\circ$ \EEE  with  John constant
  $\kappa_0\UUU \in (0,1]$ \EEE   one also has
  \begin{equation}
    \label{eq:hoelder}
    |U_t\setminus U_r|\leq c (t-r)^\mu \quad \forall 0<r\leq t \leq
    T
  \end{equation}
where $c>0$ and $\mu\in (0,1)$ depend on $n$, $\alpha_*$,
$\alpha^*$, $\kappa_0$, $d(x_0,\partial K)$, and $T$ only.
\end{proposition}

In particular, \eqref{eq:lipschitz} proves that $t\mapsto U_t$ is
Lipschitz continuous with respect to the Hausdorff distance, namely, $d_H(A,B) =
  \max\{\sup_{a\in A}\inf_{b\in B}|a-b|,\sup_{b\in B}\inf_{a\in
    A}|a-b|\}$. In fact, owing to \eqref{eq:lipschitz} one has that
  $$d_H(U_t,U_r) \leq \alpha^*|t-r| \quad \forall r, \, t>0.$$
  
  No Lipschitz continuity with respect to the
  measure of the symmetric difference $|U_t\Delta U_r|$ can be
  expected, as we prove by a counterexample in Section \ref{sec:regu}. Nonetheless,
  \eqref{eq:hoelder} implies the H\"older continuity
  $$|U_t\Delta U_r|\leq c |t-r|^\mu \quad \forall r,\, t \in [0,T].$$

Another application of  the uniform John regularity of Theorem
\ref{thm:main} is in \cite{morpho3}, where the Hamilton--Jacobi problem \eqref{eq:0} for the
generalized eikonal equation is used
to describe the quasistatic evolution of an elastic body growing by surface accretion
\cite{goriely}. This calls for a solution of a free-boundary problem,
as the time-evolving reference configuration $U_t$ of the body is a
priori unknown. In this setting, the uniform-in-time John
regularity of the evolving set $U_t$ allows to prove a uniform Korn
inequality, see \cite[Prop.~2.4]{morpho3}, which is paramount for
the analysis.

More generally, John
domains play a special role in the validity of various types of
functional
inequalities, including Poincar\'e \cite{Feng,Martio,Sun}, Korn
\cite{Acosta,Jiang}, Sobolev \cite{Hajlasz} and
Gagliardo--Nirenberg \cite{Wang} ones.
Theorem \ref{thm:main} proves that, under very general assumptions on
$H$, such \UUU inequalities uniformly \EEE hold on all sublevels \UUU of
$u$. Indeed, \EEE the uniform lower bound on the John constant 
\UUU allows \EEE to prove estimates on $U_t$ which are uniform with respect to
$t\in [0,T]$. 

\section{Proof of Theorem \ref{thm:main}}\label{sec:main}

Let us \UUU start \EEE  by specifying that $u\in C(\Rz^n)$ is a viscosity
solution to \eqref{eq:0} if $u=0$ on $K$ and for all $x\in
\Rz^n\setminus K$ and all $\varphi \in
C^1(\Rz^n)$ with $\varphi(x)-u(x)=\min(\varphi-u)$ ($\varphi(x)-u(x)=\max(\varphi-u)$) one has that $H(x,\nabla
\varphi(x))\leq 0$ ($H(x,\nabla
\varphi(x))\geq 0$, respectively), \UUU see \cite{Bardi,Barles,Lions}.
\EEE Under \UUU assumption \eqref{eq:H} \EEE
on $H$ one has that \UUU problem \eqref{eq:0} \EEE admits a unique
nonnegative viscosity solution \EEE \cite[Thm.~3.15]{Mennucci} \UUU
characterized \EEE by the formula

\begin{align}
     u(x) = \min \Bigg\{\int_0^{L_\rho}
  \sigma(\rho(s),\dot\rho(s)) \, {\d s} \::\:  &\rho \in W^{1,\infty} (0,L_\rho), \\
                &|\dot\rho|=1 \ \text{a.e.}, \ 
      \rho(0)=x, \ \rho(L_\rho)\in K \Bigg\}.  \label{eq:vis3}
\end{align}
For all $x\in \Rz^n$, the function $\sigma(x,\cdot):\Rz^n \to
[0,\infty)$ is defined (up to a sign) as the support function of the
closed and convex set $C_x=\{p\in \Rz^n \::\: H(x,p)\leq 0\}$, namely,
$$\sigma(x,v)= \sup\{-p\cdot v\::\: p \in C_x\}.$$
Note in particular that $\sigma(x,\cdot)$ is positively
$1$-homogeneous and that, owing to assumption \eqref{eq:H} on $H$, one has
\begin{equation}
  \frac{|v|}{\alpha^*} \leq \sigma(x,v) \leq  \frac{|v|}{\alpha_*} \quad \forall x,\, v
\in \Rz^n.\label{eq:sigma}
\end{equation}
The minimum in \eqref{eq:vis3} is always attained. In 
the following, given $x\in \Rz^n$ we call a minimizer $\rho$ in
\eqref{eq:vis3} {\it optimal curve} for $x$. \UUU Note that optimal curves are
not necessarily unique. \EEE

For the reader's convenience, we subdivide the proof of Theorem
\ref{thm:main} in
subsequent steps.

{\bf Step 1:} Let us preliminarily prove that \begin{equation}
  K + B_{t \alpha_*}(0)\subset U_t \subset
  K + B_{t \alpha^*}(0)\quad \forall t>0.\label{eq:nelcaso}
  \end{equation} 
To check that  the sets  $U_t$ are bounded for all $t>0$ one uses the
  bounds \eqref{eq:sigma} and the characterization in \eqref{eq:vis3} to get that
\begin{equation}
  \label{eq:be}
 {} \frac{ d(x,K)}{\alpha^*} \leq u(x) \leq  \frac{
   d(x,K)}{\alpha_*} \quad \forall\ x \in \Rz^n.
\end{equation}
 Indeed,  for any given $\ove x \in K$,  by considering the straight curve $$\ove \rho(t) = x
-t(x-\ove x)/|x-\ove x|$$ for $t \in [0,|x-\ove x|]$ we get
\begin{equation}
  \label{eq:A1}
   u(x)  \leq  \int_0^{|x-\ove x|}  \sigma(\ove \rho(s),
   \dot{\ove\rho}(s))\,{\rm d}s \stackrel{\eqref{eq:sigma}}{\leq}  \frac{|x-\ove x|}{\alpha_*}.
\end{equation}
 Passing to the minimum with respect to $\ove x \in K$ (which is
possible, since $K$ is compact) one gets the second inequality in
\eqref{eq:be}. 
On the other hand, letting $\rho$ be the optimal curve for $x$ we
readily get that 
\begin{equation}
  \label{eq:A2}
   u(x)  =\int_0^{L_{\rho}}\sigma(\rho(s),\dot\rho(s))\,{\rm
    d}s\stackrel{\eqref{eq:sigma}}{\geq} \frac{L_{ \rho}}{\alpha^*} \geq \frac{d(x,K)}{\alpha^*}
\end{equation}
 whence the first inequality in  \eqref{eq:be} follows.

 Inequalities  \eqref{eq:be} imply \eqref{eq:nelcaso}. 
   Indeed, if $x\in K+B_{t\alpha_*}(0)$, one uses the second
  inequality in \eqref{eq:be} to get $u(x) \leq d(x,K)/\alpha_* < t
  \alpha_*/\alpha_*=t$,   which implies that $x \in U_t$. On the
  contrary, if $x \in U_t$ the first inequality in \eqref{eq:be}
  gives $ d(x,K)  \leq \alpha^* u(x) <t \alpha^*$, namely,
  $x\in K+B_{t \alpha^*}(0)$.  Note that in case $\alpha \equiv
  1$, \UUU for the eikonal equation $|\nabla u|=1$ \EEE  inclusions
  \eqref{eq:nelcaso} \UUU imply that \EEE
  $U_t = K + B_t(0)$.\\
  
\begin{figure}
  \pgfdeclareimage[width=120mm]{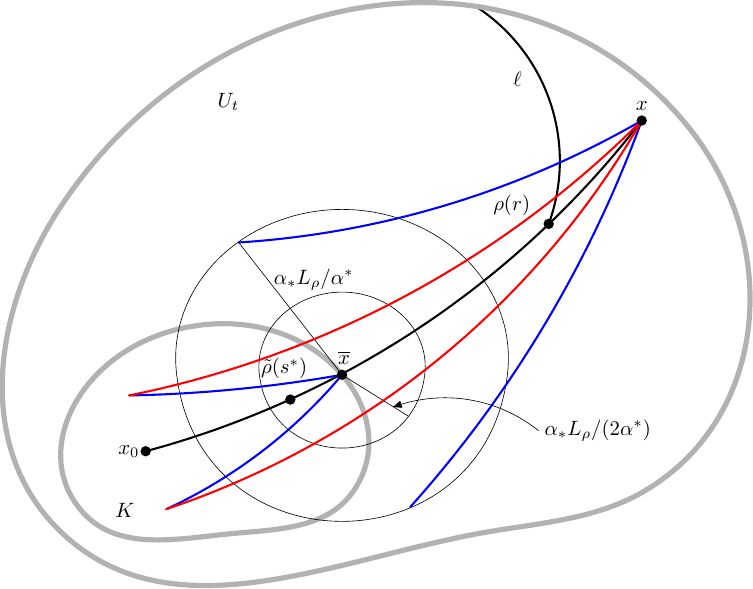}{figure_proof}
   \begin{center}
    \pgfuseimage{figure_proof}
   \end{center} 
  \caption{The construction for the proof of Theorem \ref{thm:main}}\label{figure:john}. 
\end{figure}

  \noindent\textbf{Step 2:}
We first consider the special case of $K$ being a John domain
of John constant $\kappa_0$ with respect to $x_0$ and prove that, for
all $t>0$, the sublevel $U_t$ is a John domain with respect to
$x_0$ with John
constant $\kappa$ from \eqref{eq:CJohn}.  
This follows from a geometric construction, which is illustrated in
Figure \ref{figure:john}.

Let $x \in U_t$ be given 
 and let $\rho$ be optimal for $x$, so that $\rho(L_\rho)=\ove x
\in K$.   
 The curve $\rho$ is entirely contained in $U_t$.
 More precisely, for any
$r\in [0,L_\rho]$, we  claim that 
\begin{equation}
  \label{eq:L}
  d(  \rho(r) , \partial U_t) >
  \frac{\alpha_*}{\alpha^*}r \quad \forall r \in [0,L_{\rho}].
\end{equation}

 In order to prove \eqref{eq:L},   let
$\ell:[0,L_\ell]\to \Rz^{ n}$ be any arc-length parametrized curve connecting $
\rho(r)$ to $\partial U_t$, namely, such that  $\ell(0) \in \partial U_t$ and $\ell(L_\ell)=\rho(r)$. The curve resulting from following
$\ell(s)$ for $s\in[0,L_\ell]$ and then 
$\rho(s-L_\ell+r)$ for $s \in  (  L_\ell,L_\rho+L_\ell-r]$ connects $\partial
U_t$ to  $K$.  On the one hand,  formula \eqref{eq:vis3}  and the fact that $u
= t$  on $\partial U_t$  give 
\begin{align}
  &
t \leq \int_0^{L_\ell}  \sigma(\ell(s),\dot\ell(s)) \,{\d s}+
\int_{L_\ell}^{L_\rho+L_\ell-r}  \sigma(\rho(s-L_\ell+r), \dot\rho(s-L_\ell+r))
    \,{\d s}   \nonumber\\
  &\quad =\int_0^{L_\ell}  \sigma(\ell(s),\dot\ell(s))  \,{\d s}+ \int_r^{L_\rho} \sigma(\rho(s),\dot\rho(s))  \,{\d s}
.\nonumber
\end{align}
 On the other hand, 
$$      \int_0^r   \sigma(\rho(s),\dot\rho(s)) \,{\d s}+
\int_r^{L_\rho}   \sigma(\rho(s),\dot\rho(s))  \,{\d s}=
\int_0^{L_\rho}  \sigma(\rho(s),\dot\rho(s))  \,{\d s}  =
  u(x)
 < t.$$
Putting the two inequalities together we obtain
$$ \int_0^r   \sigma(\rho(s),\dot\rho(s))  \,{\d s} < \int_0^{L_\ell} \sigma(\ell(s),\dot\ell(s))  \,{\d s}.$$
By using the bounds  \eqref{eq:sigma} on 
$ \sigma$,  the latter gives
$$ \frac{r}{\alpha^*} <\frac{L_\ell}{\alpha_*}.$$
This proves that the length $L_\ell$ of any curve connecting $ \rho(r) $ to $\partial
U_t$ is strictly bounded  from  below by $\alpha_* r/\alpha^*$. In particular, the lower
bound \eqref{eq:L} follows by considering a straight curve  from $\rho(r)$ to $\partial U_t$.   

Recall now that  $K$ is a John domain with respect to $x_0$
with John  constant $\kappa_0$ and denote by $ \rho_0: [0,L_{\rho_0}]
\to K$ an
arc-parametrized curve fulfilling $ \rho_0 (0)=\bar x$, $
\rho_0(L_{ \rho_0})=x_0$, as well as $d( \rho_0(s),\partial
K)
\geq \kappa_0s$ for all $s\in [0, L_{\rho_0}]$. We now
concatenate  the curves $ \rho $ and $ \rho_0$ to define  the
curve $\tilde \rho(s): [0,L_{\rho} +L_{\rho_0} ] \to \Rz^{ n}$ as
$$ \tilde \rho(s) =
\left\{
  \begin{array}{ll} 
     \rho(s)\quad &\text{for} \  s \in [0, L_{ \rho}],\\
      \rho_0(s-L_{\rho})\quad &\text{for} \  s \in  (  L_{ \rho} ,L_{ \rho} +L_{ \rho_0} ].
    \end{array}
  \right. 
$$
We aim at showing that, by choosing $ \kappa $ as in
\eqref{eq:CJohn},
one has that
\begin{equation}\label{eq:john}
  d(\tilde \rho(s),\partial
U_t)\geq \kappa s \quad  \text{for all} \ s \in [0,L_{ \rho}
+L_{ \rho_0} ].
\end{equation} 
In order to check \eqref{eq:john}, we distinguish three cases,
depending on the possible values of $s$ in the interval
$[0,L_{ \rho} + L_{  \rho_0}]$:

$\bullet$ Case $s\in [0,L_{ \rho}]$: Property \eqref{eq:john} follows from 
inequality \eqref{eq:L}  since
$$\kappa \stackrel{\UUU \kappa_0\leq 1}{\leq }   \frac{\alpha_*}{2\alpha^*+\alpha_*} <
\frac{\alpha_*}{\alpha^*}.$$

$\bullet$ Case $s\in [L_{ \rho}, s^*]$ with
$$s^*:=\min\left\{L_{ \rho} \left(1+
\frac{\alpha_*}{2\alpha^*}\right), L_{ \rho} + L_{\rho_0}\right\}.$$
As $s-L_\rho \leq s^*-L_\rho \leq \alpha_*L_\rho/(2\alpha^*)
$ we have that 
\begin{align}
   |\ove x - \tilde \rho(s)|&= |\ove x -  \rho_0(s-L_\rho)|
   =|\rho_0(0) -  \rho_0(\TTT s  -L_\rho)|\leq \TTT s - L_\rho \leq  {\alpha_*L_{ \rho }}/{(2\alpha^*)}, \nonumber
  \end{align}
    where we have used that $\rho_0$ is parametrized by arc-length. On
    the other hand, we have \UUU that \EEE $d(\ove x , \partial U_t) \geq  {\alpha_*L_{
    \rho }}/{\alpha^*} $ from \eqref{eq:L}. We can hence conclude that  $d(\tilde \rho(s) , \partial U_t) \geq  {\alpha_*L_{
    \rho }}/{(2\alpha^*)} $ and property \eqref{eq:john} follows 
as 
$$ \frac{\alpha_*L_{\rho }}{2\alpha^*} = \frac{\alpha_*}{2\alpha^*+\alpha_*}L_\rho\left( 1+\frac{\alpha_*}{2\alpha^*}\right)\stackrel{\eqref{eq:CJohn}}{\geq } \kappa\, s^*\geq   \kappa \,s.  $$

$\bullet$ Case $s\in [s^*, L_{ \rho}+L_{ \rho_0}]$. Note that this
case is solely relevant in case
\begin{equation}
  s^* =  L_{ \rho} \left(1+
  \frac{\alpha_*}{2\alpha^*}\right) < L_{ \rho} +
L_{\rho_0}.\label{sstar}
\end{equation}
We prove 
\eqref{eq:john} at $s\in [s^*, L_{ \rho}+L_{ \rho}]$ by using  that $K$  is a
John domain, namely, 
$$d(\tilde \rho(s), \partial U_t) > d(\tilde \rho(s), \partial
K)=d(\rho_0(s-L_\rho), \partial
K) \geq \kappa_0 (s- L_{ \rho}).$$
Property
\eqref{eq:john} then follows from
\begin{align}
  &\kappa_0 (s- L_{ \rho}) =  \kappa_0 \left( 1-
  \frac{L_{ \rho}}{s}\right)s \geq \kappa_0\left( 1-
    \frac{L_{ \rho}}{s^*}\right)s  
  = \kappa_0\frac{\alpha_*}{2\alpha^*+\alpha_*}s  \stackrel{\eqref{eq:CJohn}}{\geq }   \kappa\, s.
\end{align}

\UUU Let \EEE us remark that case
\eqref{sstar} occurs for small values of \UUU $t$ \EEE   only. Indeed, \eqref{sstar} corresponds to 
\begin{equation}\label{eq:noproof}
\frac{\alpha_*L_\rho}{2\alpha^*} \leq L_{\rho_0}.
\end{equation} 
The left-hand side of \eqref{eq:noproof} can be bounded  from 
below as  follows 
$$ \frac{\alpha_*L_\rho}{2\alpha^*}  \geq \frac{\alpha_* d(x,K)}{2\alpha^*} \stackrel{\eqref{eq:be}}{\geq}
 \frac{\alpha_* ^2}{2\alpha_*} u(x).$$
On the other hand, as  $K$   is a
John domain one has that $B_{\kappa_0 L_{\rho_0}}(x_0)\subset K$. \UUU
Recalling that \EEE 
$R = \sup\{r>0 \: : \: B_r(x_0)\subset K\}$, we have that $\kappa_0
L_{\rho_0} \leq R$ and  the right-hand side  of \eqref{eq:noproof} can
be bounded  \UUU from above \EEE  
by $R/\kappa_0$. Hence, \UUU if  $t$ is such that \EEE $ {\alpha_*^2}
\UUU t\EEE/(2\alpha_*)>  R/\kappa_0$, inequality \eqref{eq:noproof}
does not hold.  \UUU In particular, \EEE  $s^* > L_{ \rho}+L_{ \rho_0}$ and the case $s\in [s^*, L_{
  \rho}+L_{ \rho_0}]$ need not be considered. \UUU This implies that
$\kappa$ can be chosen as $\kappa=\alpha_*/(2\alpha_*+\alpha^*)$, i.e.,
independently of $\kappa_0$, for $t>2R/(\alpha_*\kappa_0)$. \EEE
\\

{\bf Step 3:} The argument of Step 2 is easily generalized
to the case of
\begin{equation}
  K\subset J \subset U_{t_0},\label{eq:step3}
\end{equation}
where $t_0>0$ and $J$ is a compact John domain  with respect to $x_0$
of John constant $\kappa_0$. In particular, for
all $t\geq t_0$ one proves that the sublevel $U_t$ is a John domain with respect to
$x_0$ with John
constant $\kappa$ from \eqref{eq:CJohn}.

The only
modification required in the proof of Step 2 is in the definition of the
concatenated curve $\tilde \rho$. Given an optimal curve $\rho$ for $x$ one
finds $\hat L\leq L_\rho$ with $\rho(\hat L)\in J$. Letting
$\rho_0:[0,L_{\rho_0}]\to \Rz^n $ be such that $\rho_0(0)=\rho(\hat
L)$, $\rho_0(L_{\rho_0})=x_0$, and $d(\rho_0,\partial J)\geq
\kappa_0s$ for all $s\in [0,L_{\rho_0}]$ one may use 
$$ \tilde \rho(s) =
\left\{
  \begin{array}{ll} 
     \rho(s)\quad &\text{for} \  s \in [0, \hat L],\\
      \rho_0(s-\hat L)\quad &\text{for} \  s \in  (  \hat L ,\hat L +L_{ \rho_0} ].
    \end{array}
  \right. 
  $$
  in order to check again \eqref{eq:john}. Note nonetheless that the
  proof in Step 2 is valid for all $t>0$, whereas under
  \eqref{eq:step3} one can prove the John regularity of $U_t$ for
  $t\geq t_0$ only.\\

  {\bf Step 4:} Consider now the general case of a compact and
  connected $K\subset \Rz^n$. Let $t_0>0$ be given and cover $K$ with a
  finite union $J$ of balls of radius smaller that
  $\alpha_*t_0$. In particular, $J$ is a John domain and, using the
  first inclusion in \eqref{eq:nelcaso} one has that \eqref{eq:step3}
  holds. One can hence follow Step 3 in order to prove that $U_t$ for
  $t \geq t_0$ is a John domain, as well. On the other hand, as $t_0$ is
  arbitrary, this argument ensures that $U_t$ is a John domain for all
  $t>0$.

  Note that, in this case, no
  lower bound on the John constant can be derived, as discussed in
  Section \ref{sec:sharpK} below.

  \section{ Sharpness of Theorem \ref{thm:main}}\label{sec:counter}

   In this section, we present some examples illustrating
  different aspects of the sharpness of our regularity result.  

In all of this section, we refer to the special case of the
generalized eikonal equation $H(x,p)=\alpha(x)|p|-1$  for different
choices of  \UUU $\alpha\in C^{0,1}(\Rz^n)$ \EEE with $0<\alpha_*\leq
\alpha(\cdot)\leq \alpha^*$. 
\UUU The unique nonnegative viscosity solution $u$ to \eqref{eq:0} \EEE
turns out to be the value function of the {\it minimal-time} problem \cite{Cannarsa,Lorenz} driven
by the controlled ODE system 
$\dot y  = \alpha(y)v$, where the measurable {\it control} $t\mapsto v(t)$ is such that
$|v|\leq1$ almost everywhere. Indeed, when $\alpha$ is Lipschitz continuous one has that the
latter differential problem admits a unique strong solution   for any
 such controls $v$  and any given initial datum $x$. One hence has that  
\begin{align}
\label{eq:min-time}
  u(x) &\UUU = \EEE\min  \Big\{t\geq 0\,:\, y(t)=x, \ \dot y
= \alpha(y)v \ \text{a.e. in} \ (0,t),  \
                                            v:(0,t)\to \Rz^{ n} \\
  &\qquad \qquad \qquad
    \text{measurable with} \ |v| \leq 1 \ \text{a.e.}, \ y(0)\in K\Big\}.
\end{align}
Correspondingly, $\overline{U_t}$ turns out to be the {\it reachable set}
at time $t$ starting from $K$ for the same controlled ODE
system, namely,
\begin{align}
\nonumber
 U_t &= \bigcup_{s < t } \Big\{y(s) \: : \: \dot y
= \alpha(y)v \ \text{a.e. in} \ (0,t),  \ \text{for some} \ v:(0,t)\to \Rz^{ n} \
    \\ \label{eq:Ut-alt} &\qquad \qquad \text{measurable with} \ |v| \leq 1 \ \text{a.e.} \ \text{and
                                            some}\  y(0)\in K\Big\}.
\end{align}
Note that the sets $U_t$ are generally
not smooth, regardless of the smoothness of $\alpha$. In the smooth
case one can equivalently  qualify  the evolution of $\partial
U_t$ in time at point $x$  as  driven by  the 
normal velocity $\alpha(x) \nu  (x)$, where $ \nu  (x)$ is the outward-pointing
normal to $\partial U_t$ at $x$.

In all examples below, we  indicate  points $x\in \Rz^{ n}$
 by
$x=(x',x_{ n})$, where $x'\in \Rz^{ n-1}$ and $x_{ n}\in \Rz$, 
and use $(e_1,\dots,e_{ n})$ for   the basis of $\Rz^{ n}$. The initial set $K$
is assumed to be smooth, contained in the half space $\{x_{ n} \leq 0\}$,
rotationally symmetric with respect to $\{x'=0\}$, and such that the
boundary $\partial K$ contains the $( n-1)$-dimensional ball
$\{(x',0)\::\: |x'|<R\}$ for $R>0$ given and large.

\subsection{Sharpness of the assumptions on $K$}\label{sec:sharpK}
Let us firstly mention that the assumption on connectedness of $K$
cannot be dropped from the statement of Theorem~\ref{thm:main}. Indeed, starting from a disconnected $K$, the second
inclusion in \eqref{eq:nelcaso} guarantees that $U_t$ is disconnected
for small $t$, hence, not a John domain.

On the other hand, for $t_0$ large enough $K+B_{t_0\alpha_*}(0)$ is necessarily
connected. Starting from such $t_0$ and arguing as in Step 4 above one
obtains that $U_t$ is John regular for all $t \geq t_0$.

\begin{figure}
  \pgfdeclareimage[width=65mm]{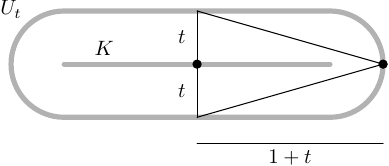}{example5}
   \begin{center}
    \pgfuseimage{example5}
   \end{center} 
  \caption{The example of Section \ref{sec:sharpK}}\label{figure:example5}
\end{figure}

For general compact and connected $K\subset \Rz^n$, although all $U_t$
are John regular, no lower bound on the John constant can be obtained. This is easily
checked by considering $K=\{x\in \Rz^n \: :\: |x_n|\leq 1\}$ and $\alpha\equiv 1 $. In this simple case,
inclusion \eqref{eq:nelcaso} entail that $U_t = K + B_t(0)$ which is
John-regular with respect to the origin with constant $t/(1+t)$ (as
well as John-regular with
respect to any other point in $U_t$ with an even smaller
constant), see Figure \ref{figure:example5}. 

\subsection{No lower bound on the interior-ball
  radius}\label{sec:counter1}
 In order to compare our result with the one in \cite{Lorenz}, 
we present an example showing that for $\alpha \not \in C^{1,1}$
one cannot expect \UUU $U_t$ to fulfill an 
interior-ball condition with some uniformly positive radius as \EEE
$t\to \infty$.  Let $0<\alpha_*< \alpha^*$, and  consider  the   Lipschitz continuous function 
\begin{equation}
\alpha(x) =
\left\{
  \begin{array}{ll}
    \alpha^* \quad&\text{ for} \ \ |x'|\leq \delta,\\[1mm]
    \alpha^*\left(1-
    \displaystyle\frac{|x'|-\delta}{\epsi\delta}\right)+\alpha_* \displaystyle\frac{|x'|-\delta}{\epsi\delta}&\text{ for}
                                                                 \ \
                                                                 \delta
                                                                 <|x'|< (1+\epsi)\delta,\\[1mm]
    \alpha_*\quad&\text{ for} \ \ |x'| \geq 
                   (1+\epsi)\delta. 
  \end{array}
\right.\label{eq:ex1}
\end{equation}
 The parameter  $\delta>0$ will be chosen later (in relation with
$t>0$),  while 
$\epsi>0$ is arbitrarily small. Note that $\alpha$ is Lipschitz
continuous with $\|\nabla \alpha\|_\infty = (\alpha^* -
\alpha_*)/(\epsi\delta)$.

 In what follows, we assume that
\begin{equation}
   \alpha^* >\sqrt{5}\alpha_*\label{eq:condition}
  \end{equation}
as this simplifies computations. Note however that
\eqref{eq:condition} could  be removed and one could treat the general case
$\alpha^*>\alpha_*$ as well, at the expense of a more involved
argument.

\begin{figure}
  \pgfdeclareimage[width=90mm]{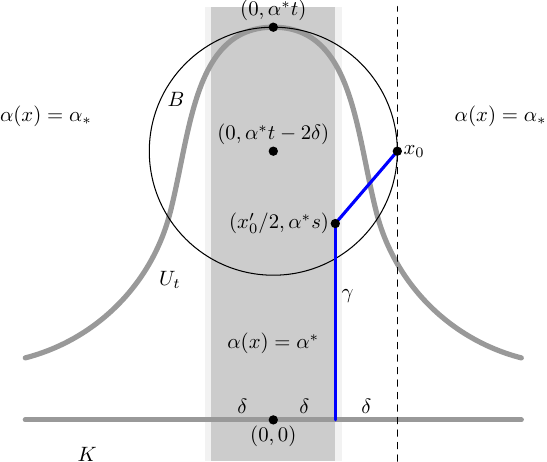}{example1}
   \begin{center}
    \pgfuseimage{example1}
   \end{center} 
  \caption{The example of Section \ref{sec:counter1}}\label{figure:example1}
\end{figure}

As $\alpha$ is Lipschitz continuous  and
radially symmetric with respect to $x'$, the \UUU sets \EEE $U_t$ are
well-defined and
radially symmetric, as well. In particular, due to radial symmetry  and \eqref{eq:vis3}  one has
that $u(0,\alpha^*t)=t$ for every $t \geq 0$, namely, that  the point
$(0,\alpha^*t)$ \UUU belongs to \EEE $ \partial U_t$ for all $t \UUU >
\EEE 0$,  recall that $K\subset U_t$ and 
see Figure \ref{figure:example1}.  In fact,  the only optimal
curve for $(0,\alpha^*t)$ is a straight segment to $(0,0)$.  

We  first show  that, for all $t>0$ there exists $\delta>0$ large enough, so that $U_t$ does not contain a ball with radius
$2\delta$ touching $\partial U_t$ only at $(0,\alpha^*t)$. Note that,
by 
symmetry, such ball  \UUU should  necessarily be \EEE $B:=B_{2\delta}(0,\alpha^*t
-2\delta)$.
We argue as follows: we consider a point $x_0=(x'_0,\alpha^*t - 2\delta)$ with
$|x'_0|=2\delta$, which belongs to $\ove B$, and prove that
$u(x_0) \geq t$.   \UUU This will in particular entail \EEE that   $x_0  \not \in  \overline{U_t}$,  hence $B \not
\subset U_t$.

In order to check that $u(x_0) \geq t$ one has to prove that no curve $\gamma$
can connect $K$ and  $x_0$ in time  smaller \UUU than \EEE  $t$. More
precisely,  cf. \eqref{eq:Ut-alt},  we \UUU need to \EEE show that
for all $\gamma_0 = \gamma(0)\in K$ and all measurable $v :(0,t) \to
\Rz^{ n}$ with $|v|\leq 1$ and $\dot \gamma=\alpha(\gamma) v $ a.e. one has
that $\gamma(t)\not =x_0$. Indeed, if $\gamma$
is such that $|(\gamma(\cdot))'| \geq  (1+\epsi)\delta $  for all times (recall that
the notation $\gamma'$ refers to the first $ n -1$ components of the
vector $\gamma$),  then, in view of \eqref{eq:min-time}  the minimal
time to connect $x_0$ with $K$ (recall that  $\{(x',0)\: : \:
|x'|<R\}\subset \partial K$ with $R$  large) is $
(\alpha^* t - 2\delta)/\alpha_*$. This is however strictly larger than 
$t$ for
\begin{equation}\label{eq:deltat}
  (\alpha^*-\alpha_*)t >2\delta ,
\end{equation}
that is, for $\delta $ large enough,
given $t$. As a consequence, no curve with $|(\gamma(\cdot))'| \geq
 (1+\epsi) \delta $  for all times can connect  $x_0 $ to
 $K$.

We are left with the case of a curve $\gamma$ with $|\gamma(\cdot)'| <
 (1+ \epsi)\delta$  at least for some
times. As $\epsi$ is assumed to be arbitrarily small, 
in the simple geometry of this example the  optimal curve for $x_0$ is
arbitrarily close \UUU (with respect to $\epsi \to 0$) \EEE to  the curve  $\gamma$ following the line 
$\{x'=x_0'/2\}$ up to some time $s<t$ and then the   straight segment 
between $(x_0'/2,\alpha^*s)$ and $x_0$,  namely, 
$$\gamma( \tau)=
\left\{
  \begin{array}{ll}
    (x_0'/2,\alpha^* \tau)\quad&\text{for} \ \  0\leq \tau <s,\\[3mm]
   \displaystyle \frac{t-\tau}{t-s}  (x_0'/2,\alpha^* s) + \left( 1-
    \displaystyle\frac{t-\tau}{t-s}\right) x_0 \quad&\text{for} \ \  s\leq \tau \leq t.
  \end{array}
\right.
$$ 
Note that the distance between
$ \gamma(s)= (x_0'/2,\alpha^*s)$ and $x_0$ is
$((\alpha^*(t-s)-2\delta)^2 +\delta^2)^{1/2}$.
In order to cover such segment  in
time $t-s$ one would need 
\begin{align}
  &\Big((\alpha^*(t-s)-2\delta)^2 +\delta^2\Big)^{1/2}  =
    \int_0^{t-s} \alpha(\gamma(\tau))\, \d \tau \nonumber\\
  &\quad  = \int_0^{t-s}\max\left\{-\frac{\alpha^*- \alpha_*}{\epsi
    (t-s)}\tau+\alpha^*, \alpha_* \right\}\, \d \tau   < (\epsi\alpha^*
  + (1{-}\epsi)\alpha_*) (t-s).
\end{align} 
\TTT The integrand featuring the $\max$ above corresponds to the velocity at time $\tau$. This is the constant $\alpha_*$ for $\tau>\epsi (t-s)$ and a convex combination between $\alpha^*$ and $\alpha_*$ for $\tau \in [0,\epsi (t-s)]$.
The  inequality between the first and the last term above  is equivalent to
\begin{equation}\label{eq:fate}
  \left((\alpha^*)^2 - (\epsi \alpha^* + (1{-}\epsi)\alpha_*)^2\right)(t-s)^2 + 5\delta^2 < 4\delta
  \alpha^*(t-s).
  \end{equation}
\UUU However, this cannot hold, \EEE regardless of the value $s<t$, as soon as $\epsi$ is
chosen  small enough. Indeed, we can equivalently rewrite \UUU
\eqref{eq:fate} \EEE  as
$$\left( (\alpha^*)^2 - (\epsi \alpha^* + (1-\epsi)\alpha_*)^2
 -  \frac45  {(\alpha^*)^2} \right)(t-s)^2 +
\left( \frac{2 }{\sqrt{5}} \alpha^* (t-s) - \sqrt{5}\delta \right)^2<0.$$
As the second term in the left-hand side above is nonnegative, \UUU 
inequality \eqref{eq:fate}  \EEE is false as soon as the coefficient of $(t-s)^2$ in the
first term is positive, namely,
$$(\alpha^*)^2 - (\epsi \alpha^* + (1-\epsi)\alpha_*)^2 -
 \frac{4}{5}  (\alpha^*)^2>0.$$
\UUU This however \EEE follows from \eqref{eq:condition} for $\epsi$ small
enough.

We have eventually  proved that $x_0\in \partial B$ cannot be reached
from  $K$ in
time $t$, which entails that $x_0\not \in U_t$ and, ultimately, $B \not \subset U_t$. In particular,
$U_t$ does not fulfill the interior-ball condition with radius
$2\delta$.

Assume now on the contrary that $\delta$ is given. \UUU By arguing
from \EEE \eqref{eq:deltat} \UUU we deduce \EEE that $U_t$ does not fulfill the interior-ball condition with radius
$2\delta$ for $t>0$ large enough. In particular, this entails that no
 uniform-in-time  lower bound on the  interior-ball radius
 of $U_t$ can be inferred \UUU in case $\alpha \in
 C^{0,1}(\Rz^n)$. \EEE 

\subsection{No interior-ball regularity}\label{sec:counter3}

Let us now consider problem \eqref{eq:0} in the open set
$\Rz^{ n-1}\times \{x_{ n}<\alpha^*T\}$ for some  arbitrary but 
  fixed final time $T>0$. For any
$\alpha^* > \alpha_*>0$, \UUU we \EEE  choose
$0<\beta<\eta<\pi/2$ in such a way that
\begin{equation}
  \label{eq:beta}
  \frac{\sin(\eta-\beta)}{\cos\beta} \geq \frac{\alpha_*}{\alpha^*}
\end{equation}
by possibly taking $\eta$ large and $\beta$ small enough. Moreover,
\UUU we \EEE
let $\epsi>0$ be arbitrarily small with $(1{+}\epsi)\beta <\eta$.

\begin{figure}
  \pgfdeclareimage[width=105mm]{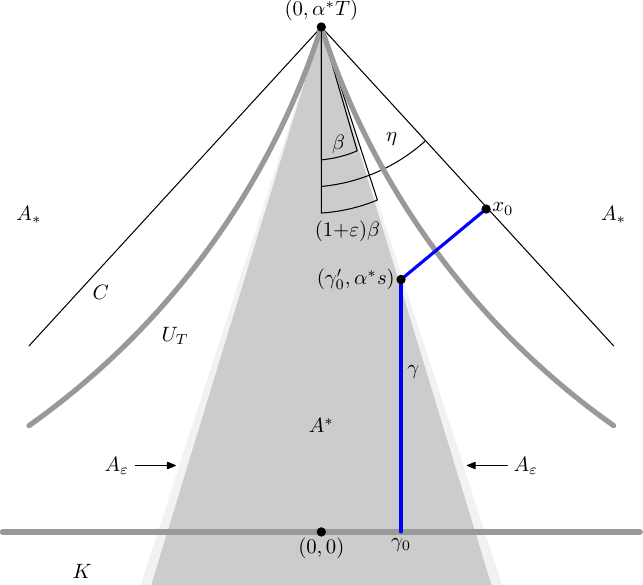}{example3}
   \begin{center}
    \pgfuseimage{example3}
   \end{center} 
  \caption{The example of Section \ref{sec:counter3}}\label{figure:example3}
\end{figure}

\UUU We define \EEE  the
Lipschitz continuous   function 
\begin{equation}
\alpha(x) =
\left\{
  \begin{array}{ll}
    \alpha^* \quad\text{if} \ x\in  A^*, \\[1mm]
     \left(1-\theta(x)\right)\alpha^*+\theta(x)\alpha_*  &\text{if}  \
              x\in A_\epsi,\\[1mm]
    \alpha_*\quad\text{if} \ x\in  A_*,
  \end{array}
\right.\label{eq:ex3}
\end{equation}
where the  regions $A^*$, $A_\epsi$,  and $A_*$ are \UUU given by \EEE
\begin{align}
  &A^*:=\left\{(x',x_{ n})\in \Rz^{ n} \times \{x_n
    <\alpha^*T\}  \: : \:\frac{|x'|}{\alpha^*T-x_{ n}}\leq \tan\beta
    \right\},\\
  &A_\epsi:=\left\{(x',x_{ n})\in \Rz^{ n} \times \{x_n
    <\alpha^*T\}   \: : \:\tan \beta
    <\frac{|x'|}{\alpha^*T-x_{ n}} < \tan((1{+}\epsi)\beta)
    \right\},\\
  & A_*:=\left\{(x',x_{ n})\in \Rz^{ n} \times \{x_n
    <\alpha^*T\}   \: : \: \tan((1{+}\epsi)\beta)\leq
    \frac{|x'|}{\alpha^*T-x_{ n}} 
    \right\},  
\end{align}
and $\theta(x)\in (0,1) $ for $x\in A_\epsi$ is given by
$$\theta(x) = \displaystyle\frac{\displaystyle\frac{|x'|}{\alpha^*T-x_{ n}}-\tan\beta}{\tan((1{+}\epsi)\beta)-\tan \beta},$$
see Figure \ref{figure:example3}. 

By symmetry,  and by \eqref{eq:vis3}  one has
that $u(0,\alpha^*T)=T$, namely that  $(0, \alpha^*T) \in \partial U_T$.  Additionally, we have that $(0,\alpha^*t)\in U_T$ for every $0\leq t<T$.  We 
aim at showing  that $U_T \subset C$ locally, where $C$
is the cone with vertex at  $(0, \alpha^*T)$, axis 
 $-e_{ n}$, and opening $\eta$, namely,
$$C=\left\{(x',x_{ n})\in \Rz^{ n} \times \{x_n
    <\alpha^*T\}  \: : \:
    \displaystyle\frac{|x'|}{\alpha^*T-x_{ n}} < \tan\eta\right\},$$
see Figure \ref{figure:example3}.
In order to prove that $U_T \subset C$ locally, we show that all $x_0\in
\partial C$ with $x_0\not = (0, \alpha^*T)$ and $(x_{0})_{ n}>\alpha_*T$ cannot be reached in time $T$
from $K$. 

Consider a curve $\gamma$ connecting $\gamma_0:=\gamma(0)\in K$ with $x_0$. Assume first that
$\gamma(t)$ is in region  $A_*$ for all times. Then necessarily
$(\gamma(t))_{ n}\leq\alpha_*t$ for all times $t\in [0,T]$. In particular
$(\gamma(T))_{ n}\leq\alpha_*T < (x_0)_{ n}$ so that $x_0$ cannot be reached
by $\gamma$ in time $T$.

Assume on the contrary that $\gamma$ enters the region $A^*\cup
A_\epsi$ for some times. As $\epsi$ is arbitrarily small, in the geometry
of this example, one has that the optimal curve in this class is
arbitrary close to  the curve  $\gamma$ following the line
$\{x'=\gamma_0'\}$ with $\gamma_0'=x_0'\alpha^*(T-s)\sin\beta/|x_0'|$
up to some time $s<T$ and then  the   straight segment
between $(\gamma_0',\alpha^*s)$ and $x_0$,  namely,
$$
\gamma(t) =
\left\{
  \begin{array}{ll}
    \gamma_0 + \alpha^*t e_n\quad&\text{for} \ \ 0\leq t <s,\\[3mm]
   \displaystyle \left( \frac{T-t}{T-s}\right) (\gamma_0',\alpha^*s) +
    \left(1- \frac{T-t}{T-s}\right) x_0 \quad &\text{for} \ \ s\leq t
                                                \leq T.
  \end{array}
  \right.
  $$ 
  The distance between
$ \gamma(s) =  (\gamma_0',\alpha^*s)$ and $x_0$ is  
$$ \alpha^* (T-s) \frac{\sin (\eta-\beta)}{\cos
  \beta}.$$
On the other hand, by using condition \eqref{eq:beta} and choosing
$\epsi$ small enough depending on $\alpha_*$, $\alpha^*$, $\eta$, and $\beta$ one finds that
the latter quantity is strictly greater than
$$
\left(\displaystyle\frac{\sin(\eta-(1{+}\epsi)\beta)}{\sin(\eta-\beta)}
  \alpha_*
+ \left(1-
  \displaystyle\frac{\sin(\eta-(1{+}\epsi)\beta)}{\sin(\eta-\beta)}\right)
\alpha^* \right)(T-s),$$
which
bounds from above the
maximum distance that the curve  $\gamma$  can travel in region $A_*\cup A_\epsi$ in time
$T-s$.
We hence conclude that the curve $\gamma$ cannot reach $x_0$ in time
$T$. This proves that $U_T$ is locally contained in the cone $C$. In
particular, $U_T$ is not interior-ball regular.

Before closing this section, let us remark that $\alpha$ is defined in
$\Rz^{ n-1}\times  \{x_n <\alpha^*T \}$  only and cannot be continuously extended to the whole
$\Rz^{ n}$.

\subsection{No interior-cone regularity}\label{sec:counter4}

The example of Section \ref{sec:counter3} is not interior-ball
regular but still interior-cone regular.  \UUU In this section, we
\EEE  rework the example in
order to prove  that  even interior-cone regularity may fail.  To this aim, we
still use the choice of $\alpha,\,\beta,\,\eta$, and $\theta$  from \eqref{eq:ex3}, but redefining the
regions $A^*$, $A_\epsi$, and $A_*$  as in Figure
\ref{figure:example4} (see details below). In this case, one can prove
that  
there exists a time $T^*$ such that $U_{T^*}$ is not
interior-cone regular, but
merely John regular.

\pgfdeclareimage[height=50mm]{V}{V2}
 \pgfdeclareimage[height=50mm]{Q}{Q2}
  \begin{figure}[h]
    \begin{center}
      \pgfuseimage{V} \hspace{6mm}\pgfuseimage{Q}
    \end{center}
    \caption{A two-dimensional representation of the sets $V_\beta$ and $Q_\beta$}
    \label{fig:V}
  \end{figure} 

In order to give some  details,  let us start by defining the
 union of balls 
\begin{align}
  V_\beta:= \bigcup_{s\in[0,1]} B_{(\tan \beta) (2-s)}(se_n) \label{eq:Vbeta}
\end{align} 
see Figure \ref{fig:V}.
The parameter  $\beta>0$ plays the role of the opening of a cone, see Section
\ref{sec:counter3}, and will later be chosen to be
small. Correspondingly, for $\beta<\eta<\pi/2$ and $\epsi>0$ arbitrarily small we analogously
define $V_\eta$ and $V_{(1{+}\epsi)\beta}$.

Moreover, we let
$$Q_\beta:=V_\beta \cup  \left(\frac12   J   V_\beta  + (0,
  \dots, 0, 1 ))\right),$$ where
$ J   $ is a rotation, mapping $e_{ n}$ to $e_1$,  see
Figure \ref{fig:V}.  Correspondingly, for $\beta<\eta<\pi/2$ and $\epsi>0$ arbitrarily small we analogously
define $Q_\eta$ and
$Q_{(1{+}\epsi)\beta}$ by replacing $\beta$ by $\eta$ or
$(1{+}\epsi)\beta$, respectively.

Finally, we let $A^*$ and $A_\epsi$ be
given by
\begin{align}
  &A^* := \bigcup_{i =0}^\infty \left(\frac{1}{4^{i}}Q_\beta
+\sum_{j=1}^i\frac{1}{4^{j-1}} \left(\frac{1}{2},0,\dots,0,1\right)
    \right),\\
   &A_\epsi := \left(\bigcup_{i =0}^\infty \left(\frac{1}{4^{i}}Q_{(1{+}\epsi)\beta}
+\sum_{j=1}^i\frac{1}{4^{j-1}} \left(\frac{1}{2},0,\dots,0,1\right)
    \right)\right)\setminus A^*.\\
\end{align}
 One has that $A^*, \, A_\epsi\subset \Rz^{n-1} \times \{x_n <
  4/3\}$ and defines $A_* =\Rz^{n-1} \times \{x_n <
 \UUU 4\alpha^*/3\EEE\} \setminus (A^*\cup A_\epsi)$.

Consider now the point 
$$x^*=     \sum_{j=1}^\infty \frac{1}{4^{j-1}}\left(\frac{1}{2},0,\dots,0,1
\right) = \left(\frac{2}{3},0,\dots,0,\frac{  4 }{3}\right).$$
 Clearly, $x^* \in \overline{A^*}$ can be reached from $K$ in finite time $T^*$ by a curve entirely
in $\overline{A^*}$, see Figure \ref{figure:example4}. Assume now that $\alpha^*>\alpha_*>0$ are given in such a
way that one can choose $\eta$ with $\eta>\beta>0$ small and still
fulfilling \eqref{eq:beta} (this is,  for instance,   the case for $\alpha_* <<
\alpha^*$)  and define
$$
\alpha(x)=
\left\{
  \begin{array}{ll}
     \alpha_*\quad&\text{for} \ \ x\in A_*,\\[2mm]
    \displaystyle\frac{ \alpha^* \, d(x,A_*) +  \alpha_*  \, d(x,A^*)}{\max\{d(x,A_*)
    , d(x,A^*)\}}&\text{for} \ \ x \in A_\epsi,\\[4mm]
    \alpha^*\quad&\text{for} \ \ x\in A^*.\\
  \end{array}
\right.
$$
Note that $\alpha$ is locally Lipschitz continuous in $\Rz^{n-1}\times
\{x_n <  4/3 \}$ as one has that
$\max\{d(x,A_*)
    , d(x,A^*)\}$ is uniformly positive in $A_\epsi$. 

\begin{figure}
  \pgfdeclareimage[width=77mm]{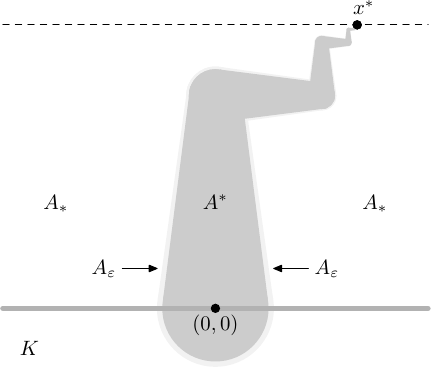}{example4}
   \begin{center}
    \pgfuseimage{example4}
   \end{center} 
  \caption{The example of Section \ref{sec:counter4}  in $\Rz^2$}\label{figure:example4}
\end{figure}

 For $\epsi$ small one can
then argue as in Section \ref{sec:counter3} in order to check that  in a
neighborhood of $x^*$  the
superlevel $U_{T^*}$ is contained in the set
$$ A_\eta:= \bigcup_{i =0}^\infty \left(\frac{1}{4^{i}}Q_\eta
+\sum_{j=1}^i \frac{1}{4^{j-1}}\left(\frac12,0,\dots,0,1\right)
\right).$$
 Since  $A_\eta$ is (John regular but) not interior-cone regular, so is $U_{T^*}$,  for $x^*\in A_\eta\cap U_{T^*}$. 

Once again, \UUU the function \EEE $\alpha$ cannot be continuously extended to $\Rz^{ n}$.


\section{Proof of Proposition \ref{prop:regu}}\label{sec:regu}

To prove the lower bound \eqref{eq:lipschitz0}, fix $x\in \partial U_t$ and $y \in
\partial U_r$, let $\rho\in W^{1,\infty}(0,L_\rho)$ be
optimal  for $y$, and define the curve $\tilde \rho\in W^{1,\infty}(0,L_\rho+ |x-y|)$ as
$$
\tilde \rho(s) =
\left\{
  \begin{array}{ll}
  x+ s\displaystyle\frac{x-y}{|x-y|} \quad &\text{for} \ s \in [0,|x-y|],\\[2mm]
  \rho(s- |x-y|) &\text{for} \ s \in (|x-y|,L_\rho +|x-y|].
  \end{array}
  \right.
  $$
  As $\tilde \rho$ connects $x$ and $K$ we have that
  \begin{align}
    & t=u(x)\leq \int_0^{L_\rho +|x-y|}\sigma(\tilde
    \rho(s),\dot{\tilde \rho}(s))\, \d s \nonumber\\
    &\quad = \int_0^{L_\rho}\sigma(
    \rho(s),\dot\rho(s))\, \d s + \int_{0}^{|x-y|}\sigma(\tilde
    \rho(s),\dot{\tilde \rho}(s))\, \d s \nonumber\\
    &\quad = u(y) + \int_{0}^{|x-y|}\sigma(\tilde
      \rho(s),\dot{\tilde \rho}(s))\, \d s\stackrel{\eqref{eq:sigma}}{\leq}  r  + \frac{|x-y|}{\alpha_*}\nonumber
  \end{align}
and \eqref{eq:lipschitz0} follows by taking the infimum first on $y\in
\partial U_r$ and then on $x\in \partial U_t$.

Let us now check the Lipschitz bound \eqref{eq:lipschitz}. Assume \UUU
with no loss of generality \EEE that 
$r<t$ and fix $x \in
U_t\setminus U_r$. Let $\rho\in W^{1,\infty}(0,L_\rho)$ be optimal for
$x$ and $\ell \in (0,L_\rho)$ be such that $r-\epsi< u(\rho(\ell))<r$
for $\epsi>0$ small. Such a
value $\ell$ exists as $ u$ is continuous along the curve $\rho$. Then,
we have
\begin{align}
  &d(x,U_r)\leq |x-\rho(\ell)| \leq L_\rho - \ell\stackrel{\eqref{eq:sigma}}{\leq} 
 \alpha^*\int_\ell^{L_\rho}\sigma(\rho(s),\dot \rho(s))\,
  \d s \nonumber\\
  &\quad =\alpha^*\int_0^{L_\rho}\sigma(\rho(s),\dot \rho(s))\,
  \d s-\alpha^*\int_0^{\ell}\sigma(\rho(s),\dot \rho(s))\,
    \d s  \nonumber\\[2mm]
  &\quad \leq \alpha^*( u(x) -
     u(\rho(\ell))) \leq \alpha^*(t - r+\epsi).\nonumber
\end{align}
By taking the infimum for $\epsi>0$ and then the supremum on $x\in
U_t$ we obtain \eqref{eq:lipschitz}.

Note that no Lipschitz continuity with respect to the measure of the
difference \UUU $U_t\setminus U_r$ \EEE of the form $|U_t\setminus
\UUU U_r\RRR |\leq \lambda |t-\UUU r \EEE |$ for some $\UUU \lambda \EEE>0$
can be
expected in general. To see this, consider the case of the generalized
eikonal equation and assume to be given $K$ and $\gamma$ in such a way
that $\partial U_t$ is Lipschitz for all $t \in (0,T)$ and
$\mathcal{H}^{n-1}(\partial U_t)\geq c/(T-t)$ for some $c>0$. This  
can be easily realized by arguing as in Section \ref{sec:counter3} for $\partial U_T$ having ${\rm dim}_M
(\partial U_T)>n-1$. As $\partial U_t$ is Lipschitz one has
\begin{equation}
  \lim_{\epsi\to 0}\frac{|(U_t + B_\epsi(0))\setminus U_t|}{\epsi} =
\mathcal{H}^{n-1}(\partial U_t)\geq \frac{c}{T-t} \quad \forall t\in
(0,T).\label{eq:nolip}
\end{equation}
On the other hand, also using \eqref{eq:lipschitz0} for all $0<r<t$ we have 
$$\liminf_{t\searrow r}\,|U_t\setminus U_r|\geq \liminf_{t\searrow
  r}\,|(U_r+B_{\alpha_*(t-r)}(0))\setminus U_r|
\stackrel{\eqref{eq:nolip}}{\geq} \liminf_{t\searrow
  r}\frac{c\alpha_*}{2(T-r)}(t-r).$$
This proves that, in order  $|U_t\setminus U_r|\leq \lambda |t-r|$ to
hold for all $r,\,t\in (0,T)$ one has to require that $\lambda \geq
\limsup_{t\in (0,T)}c\alpha_*/(2(T-t))=\infty$, \UUU which contradicts
$\lambda \in \Rz_+$. \EEE

In order to prove the H\"older bound \eqref{eq:hoelder} we use the
following Lemma, which follows directly from the proof of \cite[Cor.~2]{Smith}.

\begin{lemma}
Let $U\subset \Rz^n$ be a John domain with respect to $x_0 \in {\rm
  int} \,U$ with John constant $\kappa_0>0$. \UUU There exist
$\epsi_0>0$,  $\mu\in (0,1)$ and $\haz c>0$ depending on $n$, $d(x,\partial U)$,
and $\kappa_0$ only, such that 
the minimal number $N(\epsi)$ of balls with radius $\epsi\in (0,\epsi_0]$ needed to
cover $\partial U$ can be  bounded  by $ N(\epsi)\leq \haz c
\epsi^{\mu-n}$. \EEE
\end{lemma}

\UUU Let $0\leq r< t \leq T$ and use \EEE the lemma to find a finite covering $A$ of $\partial U_t$ by
\UUU $N(\delta)$ balls of radius
$$\delta:=\epsi_0\frac{\alpha^*(t-r)}{\alpha^* T}\in (0,\epsi_0].$$ \EEE
Owing to
\eqref{eq:lipschitz} we have that 
$U_t\setminus U_r \subset A$. Hence, 
\begin{align}
  |U_t\setminus U_r| &\leq |A|\leq  N(\UUU \delta)\UUU |B_1(0)|\delta^n
  \leq  \UUU\haz c  \delta^{\mu-n}|B_1(0)|\delta^n \nonumber\\
                     &\UUU = \haz c
                       |B_1(0)|\delta^\mu = \frac{\haz c |B_1(0)|
                       \epsi_0^\mu(\alpha^*)^\mu}{(\alpha^*)^\mu T^\mu}  (t-r)^\mu\nonumber
 \end{align}
  where $\mu\in (0,1)$ and $\haz c$ depend on $n$, $d(x,\partial
  U_t)$, and the John constant of $U_t$. Owing to Theorem
  \ref{thm:main} these are however bounded in terms of $d(x,K)$,
  $\alpha_*$, $\alpha^*$, and $T$, independently of $t\in [0,T]$,
  whence \eqref{eq:hoelder} follows.

\section*{Acknowledgements}
Support from the Austrian Science Fund (FWF) through
projects 10.55776/F65, 10.55776/I4354, 10.55776/I5149,
10.55776/P32788, 10.55776/I4052, 10.55776/V662, 10.55776/P35359, 10.55776/F100800, and 10.55776/Y1292
as well as from BMBWF
through the OeAD-WTZ project CZ 09/2023 is gratefully
acknowledged. \UUU For open access purposes, the authors have applied a CC BY public
copyright license to any author accepted manuscript version arising
from this submission. \EEE

\end{document}